\documentclass[11pt]{article}
\usepackage{bbm}
\usepackage{mathrsfs}
\usepackage{amscd}
\usepackage{amsmath,amsfonts,amssymb,amscd}
\usepackage{indentfirst,graphics,epsfig,psfrag}
\input{epsf}
\usepackage{ifpdf}
\usepackage{enumerate}
\usepackage{appendix}
\usepackage{enumerate}

\usepackage{lineno}

\makeatletter \@addtoreset{figure}{section} \makeatother
\makeatletter
\long\def\@makecaption#1#2{%
   \vskip 10\p@
   \setbox\@tempboxa\hbox{{#1}\ \ #2}%
   \ifdim \wd\@tempboxa >\hsize

       {#1}\ \ #2\par
   \else
       \hbox to\hsize{\hfil\box\@tempboxa\hfil}%
   \fi}
\makeatother

\newtheorem{thm}{Theorem}

\newtheorem{lem}{Lemma}

\newtheorem{obs}{Observation}
\newtheorem{pro}{Proposition}

\newcommand{\qed}{{\hfill\rule{3pt}{7pt}}}

\def\qed{\hfill \rule{4pt}{7pt}}

\begin{document}
\title{\textbf{The Steiner diameter of a graph}
\footnote{Supported by the National Science Foundation of China (No.
11161037) and the Science Found of Qinghai Province (No.
2014-ZJ-907).}}
\author{
\small  Yaping Mao\footnote{E-mail: maoyaping@ymail.com}\\[0.2cm]
\small Department of Mathematics, Qinghai Normal\\
\small University, Xining, Qinghai 810008, China\\
\small }
\date{}
\maketitle
\begin{abstract}
The Steiner distance of a graph, introduced by Chartrand,
Oellermann, Tian and Zou in 1989, is a natural generalization of the
concept of classical graph distance. For a connected graph $G$ of
order at least $2$ and $S\subseteq V(G)$, the \emph{Steiner
distance} $d(S)$ among the vertices of $S$ is the minimum size among
all connected subgraphs whose vertex sets contain $S$. Let $n,k$ be
two integers with $2\leq k\leq n$. Then the \emph{Steiner
$k$-eccentricity $e_k(v)$} of a vertex $v$ of $G$ is defined by
$e_k(v)=\max \{d(S)\,|\,S\subseteq V(G), \ |S|=k, \ and \ v\in S
\}$. Furthermore, the \emph{Steiner $k$-diameter} of $G$ is
$sdiam_k(G)=\max \{e_k(v)\,|\, v\in V(G)\}$. In 2011, Chartrand,
Okamoto and Zhang showed that $k-1\leq sdiam_k(G)\leq n-1$. In this
paper, graphs with $sdiam_3(G)=2,3,n-1$ are characterized,
respectively. We also consider the Nordhaus-Gaddum-type results for
the parameter $sdiam_k(G)$. We determine sharp upper and lower
bounds of $sdiam_k(G)+sdiam_k(\overline{G})$ and $sdiam_k(G)\cdot
sdiam_k(\overline{G})$ for a graph $G$ of order $n$. Some
graph classes attaining these bounds are also given.\\[2mm]
\end{abstract}

\section{Introduction}

All graphs in this paper are undirected, finite and simple. We refer
to \cite{Bondy} for graph theoretical notation and terminology not
described here. Distance is one of the most basic concepts of
graph-theoretic subjects. If $G$ is a connected graph and $u,v\in
V(G)$, then the \emph{distance} $d(u,v)$ between $u$ and $v$ is the
length of a shortest path connecting $u$ and $v$. If $v$ is a vertex
of a connected graph $G$, then the \emph{eccentricity} $e(v)$ of $v$
is defined by $e(v)=\max\{d(u,v)\,|\,u\in V(G)\}$. Furthermore, the
\emph{radius} $rad(G)$ and \emph{diameter} $diam(G)$ of $G$ are
defined by $rad(G)=\min\{e(v)\,|\,v\in V(G)\}$ and $diam(G)=\max
\{e(v)\,|\,v\in V(G)\}$. These last two concepts are related by the
inequalities $rad(G)\leq diam(G) \leq 2 rad(G)$. The center $C(G)$
of a connected graph $G$ is the subgraph induced by the vertices $u$
of $G$ with $e(u)=rad(G)$. Recently, Goddard and Oellermann gave a
survey paper on this subject, see \cite{Goddard}.

The distance between two vertices $u$ and $v$ in a connected graph
$G$ also equals the minimum size of a connected subgraph of $G$
containing both $u$ and $v$. This observation suggests a
generalization of the classical graph distance. The Steiner distance
of a graph, introduced by Chartrand, Oellermann, Tian and Zou
\cite{Chartrand} in 1989, is a natural and nice generalization of
the concept of classical graph distance. For a graph $G(V,E)$ and a
set $S\subseteq V(G)$ of at least two vertices, \emph{an $S$-Steiner
tree} or \emph{a Steiner tree connecting $S$} (or simply, \emph{an
$S$-tree}) is a such subgraph $T(V',E')$ of $G$ that is a tree with
$S\subseteq V'$. Let $G$ be a connected graph of order at least $2$
and let $S$ be a nonempty set of vertices of $G$. Then the
\emph{Steiner distance} $d_G(S)$ among the vertices of $S$ (or
simply the distance of $S$) is the minimum size among all connected
subgraphs whose vertex sets contain $S$. When there is no
$S$-Steiner tree, we set $d_G(S)=\infty$ by convention. Note that if
$H$ is a connected subgraph of $G$ such that $S\subseteq V(H)$ and
$|E(H)|=d(S)$, then $H$ is a tree. Clearly,
$d(S)=\min\{e(T)\,|\,S\subseteq V(T)\}$, where $T$ is subtree of
$G$. Furthermore, if $S=\{u,v\}$, then $d(S)=d(u,v)$ is nothing new
but the classical distance between $u$ and $v$. Clearly, if $|S|=k$,
then $d(S)\geq k-1$. If $G$ is the graph of Figure 1 $(a)$ and
$S=\{u,v,x\}$, then $d(S)=4$. There are several trees of size $4$
containing $S$. One such tree $T$ is also shown in Figure 1 $(b)$.
This example is from \cite{Chartrand}.
\begin{figure}[!hbpt]
\begin{center}
\includegraphics[scale=0.7]{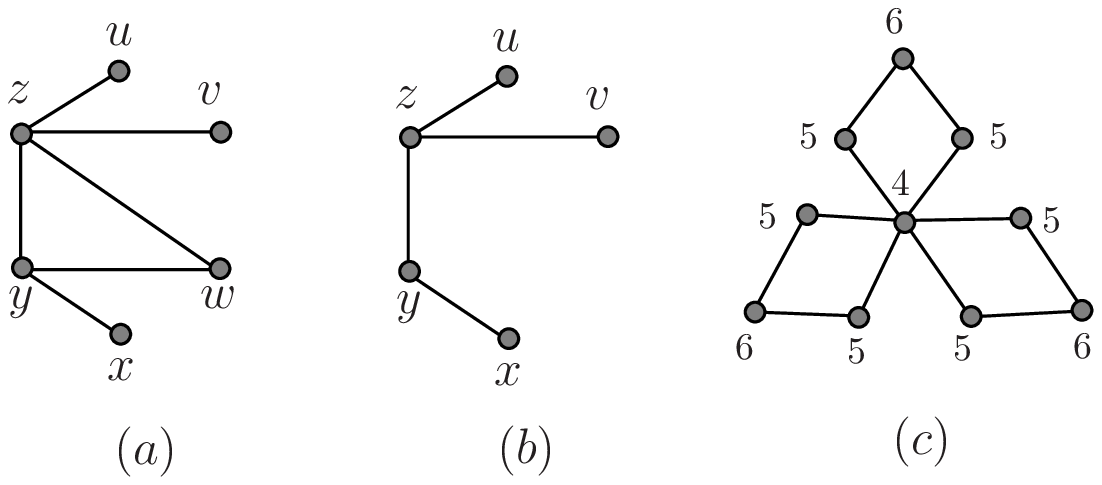}
\end{center}
\begin{center}
Figure 1: Graphs for the basic definition.
\end{center}\label{fig1}
\end{figure}

Let $n$ and $k$ be two integers with $2\leq k\leq n$. The
\emph{Steiner $k$-eccentricity $e_k(v)$} of a vertex $v$ of $G$ is
defined by $e_k(v)=\max \{d(S)\,|\,S\subseteq V(G), |S|=k,~and~v\in
S \}$. The \emph{Steiner $k$-radius} of $G$ is $srad_k(G)=\min \{
e_k(v)\,|\,v\in V(G)\}$, while the \emph{Steiner $k$-diameter} of
$G$ is $sdiam_k(G)=\max \{e_k(v)\,|\,v\in V(G)\}$. Note for every
connected graph $G$ that $e_2(v)=e(v)$ for all vertices $v$ of $G$
and that $srad_2(G)=rad(G)$ and $sdiam_2(G)=diam(G)$. Each vertex of
the graph $G$ of Figure 1 $(c)$ is labeled with its Steiner
$3$-eccentricity, so that $srad_3(G)=4$ and $sdiam_3(G)=6$.

In \cite{DankelmannSO2}, Dankelmann, Swart and Oellermann obtained
an upper bound on $sdiam_k(G)$ for a graph $G$ in terms of the order
of $G$ and the minimum degree of $G$, that is, $sdiam_n(G)\leq
\frac{3p}{\delta+1}+3n$. Recently, Ali, Dankelmann, Mukwembi
\cite{AliDM} improved the bound of $sdiam_n(G)$ and showed that
$sdiam_n(G)\leq \frac{3p}{\delta+1}+2n-5$ for all connected graphs
$G$. Moreover, they constructed graphs to show that the bounds are
asymptotically best possible.

The Steiner tree problem in networks, and particularly in graphs,
was formulated quite recently in 1971 by Hakimi (see \cite{Hakimi})
and Levi (see \cite{Levi}). In the case of an unweighted, undirected
graph, this problem consists of finding, for a subset of vertices
$S$, a minimal-size connected subgraph that contains the vertices in
$S$. The computational side of this problem has been widely studied,
and it is known that it is an NP-hard problem for general graphs
(see \cite{HwangRW}). The determination of a Steiner tree in a graph
is a discrete analogue of the well-known geometric Steiner problem:
In a Euclidean space (usually a Euclidean plane) find the shortest
possible network of line segments interconnecting a set of given
points. Steiner trees have application to multiprocessor computer
networks. For example, it may be desired to connect a certain set of
processors with a subnetwork that uses the least number of
communication links. A Steiner tree for the vertices, corresponding
to the processors that need to be connected, corresponds to such a
desired subnetwork. The problem of determining the Steiner distance
is known to be NP-hard \cite{GareyJ}.

Let $G$ be a $k$-connected graph and $u$, $v$ be any pair of
vertices of $G$. Let $P_k(u,v)$ be a family of $k$ vertex-disjoint
paths between $u$ and $v$, i.e., $P_k(u,v)=\{p_1,p_2,\cdots,p_k\}$,
where $p_1\leq p_2\leq \cdots \leq p_k$ and $p_i$ denotes the number
of edges of path $p_i$. The \emph{$k$-distance} $d_k(u,v)$ between
vertices $u$ and $v$ is the minimum $|p_k|$ among all $P_k(u,v)$ and
the \emph{$k$-diameter} $d_k(G)$ of $G$ is defined as the maximum
$k$-distance $d_k(u,v)$ over all pairs $u,v$ of vertices of $G$. The
concept of $k$-diameter emerges rather naturally when one looks at
the performance of routing algorithms. Its applications to network
routing in distributed and parallel processing are studied and
discussed by various authors including Chung \cite{Chung}, Du et al.
\cite{Du}, Hsu \cite{Hsu, Hsu2}, Meyer and Pradhan \cite{Meyer}.

In the sequel, let $K_{s,t}$, $K_{n}$, $P_n$ and $C_n$ denote the
complete bipartite graph of order $s+t$ with part sizes $s$ and $t$,
complete graph of order $n$, path of order $n$, and cycle of order
$n$, respectively. The degree of a vertex $v$ in $G$ is denoted by
$d_G(v)$. For $S\subseteq V(G)$, we denote $G-S$ the subgraph by
deleting the vertices of $S$ together with the edges incident with
them from $G$. If $S=\{v\}$, we simply write $G-v$ for $G-\{v\}$.
Let $N_G(v)$ denote the neighbors of the vertex $v$ in $G$.

From the above definitions, the following observation is easily
seen.
\begin{obs}\label{obs1}
Let $k,n$ be two integers with $2\leq k\leq n$.

$(1)$ For a complete graph $K_n$, $sdiam_k(K_n)=k-1$;

$(2)$ For a path $P_n$, $sdiam_k(P_n)=n-1$;

$(3)$ For a cycle $C_n$, $sdiam_k(C_n)=\big
\lfloor\frac{n(k-1)}{k}\big \rfloor$.
\end{obs}

In {\upshape\cite{Chartrand}}, Chartrand et al. derived the upper
and lower bounds for $sdiam_k(G)$.

\begin{pro}{\upshape\cite{Chartrand}}\label{pro1}
Let $k,n$ be two integers with $2\leq k\leq n$, and let $G$ be a
connected graph of order $n$. Then
$$
k-1\leq sdiam_k(G)\leq n-1.
$$
Moreover, the bounds are sharp.
\end{pro}

The following observation is immediate.

\begin{obs}\label{obs2}
Let $G$ be a connected graph of order $n$. Then

$(1)$ $sdiam_2(G)=1$ if and only if $G$ is a complete graph;

$(2)$ $sdiam_2(G)=n-1$ if and only if $G$ is a path of order $n$.
\end{obs}

Let $uv$ be an edge in $G$. A \emph{double-star} on $uv$ is a
maximal tree in $G$ which is the union of stars centered at $u$ or
$v$ such that each star contains the edge $uv$. Bloom \cite{Bloom}
characterized the graphs with $sdiam_2(G)=2$.

\begin{thm}{\upshape\cite{Bloom}}\label{th1}
Let $G$ be a connected graph of order $n$. Then $sdiam_2(G)=2$ if
and only if $\overline{G}$ is non-empty and $\overline{G}$ does not
contain a double star of order $n$ as its subgraph.
\end{thm}

In this paper, we focus on the case $k=3$ and characterize the
graphs with $sdiam_3(G)=2$ in Section $2$, which can be seen as an
extension of $(1)$ of Observation \ref{obs2}.

\begin{thm}\label{th2}
Let $G$ be a connected graph of order $n$. Then $sdiam_3(G)=2$ if
and only if $0\leq \Delta(\overline{G})\leq 1$ if and only if
$n-2\leq \delta(G)\leq n-1$.
\end{thm}

We now define two graph classes. A \emph{triple-star} $H_1$ is
defined as a connected graph of order $n$ obtained from a triangle
and three stars $K_{1,a}, K_{1,b}, K_{1,c}$ by identifying the
center of a star and one vertex of the triangle, where $0\leq a\leq
b\leq c$, $c\geq 1$ and $a+b+c=n-3$; see Figure 2 $(a)$. Let $H_2$
be a connected graph of order $n$ obtained from a path $P=uvw$ and
$n-3$ vertices such that for each $x\in V(H_2)- \{u,v,w\}$,
$xu,xv,xw\in E(H_2)$, or $xu,xv\in E(H_2)$ but $xw\notin E(H_2)$, or
$xv,xw\in E(H_2)$ but $xu\notin E(H_2)$, or $xu,xw\in E(H_2)$ but
$xv\notin E(H_2)$, or $xv\in E(H_2)$ but $xu,xw\notin E(H_2)$; see
Figure 2 $(c)$.
\begin{figure}[!hbpt]
\begin{center}
\includegraphics[scale=0.7]{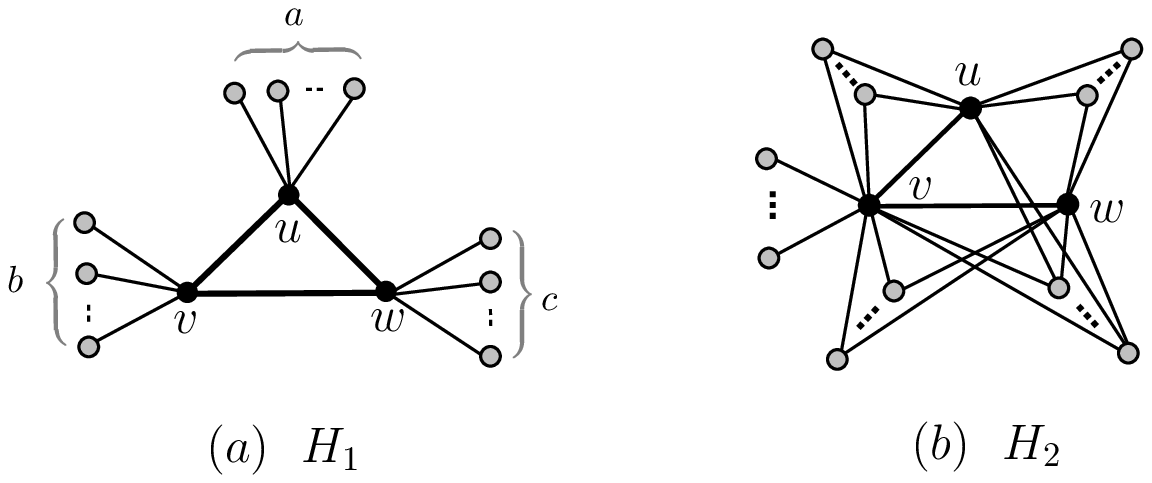}
\end{center}
\begin{center}
Figure 2: Graphs for Theorem \ref{th3}.
\end{center}\label{fig1}
\end{figure}

Graphs with $sdiam_3(G)=3$ are also characterized in Section $2$,
which can be seen as an extension of Theorem \ref{th1}.

\begin{thm}\label{th3}
Let $G$ be a connected graph of order $n$. Then $sdiam_3(G)=3$ if
and only if $G$ satisfies the following conditions.

$\bullet$ $\Delta(\overline{G})\geq 2$;

$\bullet$ $\overline{G}$ does not contain a triple-star $H_1$ as its
subgraph;

$\bullet$ $\overline{G}$ does not contain $H_2$ as its subgraph.
\end{thm}

Denote by $T_{a,b,c}$ a tree with a vertex $v$ of degree $3$ such
that $T_{a,b,c}-v=P_{a}\cup P_{b}\cup P_{c}$, where $0\leq a\leq
b\leq c$ and $1\leq b\leq c$ and $a+b+c=n-1$; see Figure 3 $(a)$.
Observe that $T_{0,b,c}$ where $b+c=n-1$ is a path of order $n$.
Denote by $\bigtriangleup_{p,q,r}$ a unicyclic graph containing a
triangle $K_{3}$ and satisfying $\bigtriangleup_{p,q,r}-
V(K_{3})=P_{p}\cup P_{q}\cup P_{r}$, where $0\leq p\leq q\leq r$ and
$p+q+r=n-3$; see Figure 3 $(b)$.
\begin{figure}[!hbpt]
\begin{center}
\includegraphics[scale=0.7]{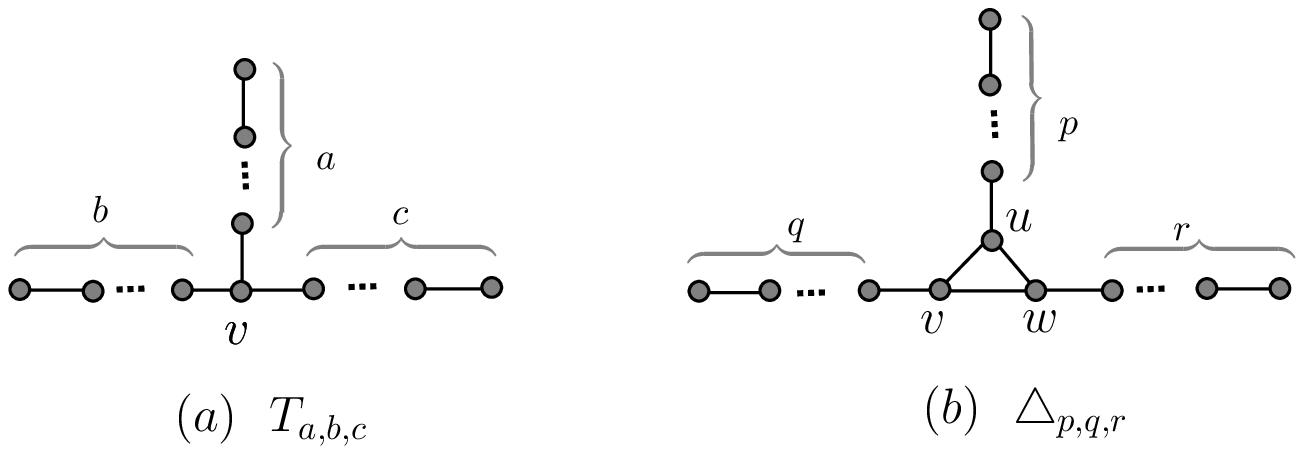}
\end{center}
\begin{center}
Figure 3: Graphs for Theorem \ref{th4}.
\end{center}\label{fig1}
\end{figure}

In Section $2$, graphs with $sdiam_3(G)=n-1$ are also characterized,
which can be seen as an extension of $(2)$ of Observation
\ref{obs2}.

\begin{thm}\label{th4}
Let $G$ be a connected graph of order $n \ (n\geq 3)$. Then
$sdiam_3(G)=n-1$ if and only if $G=T_{a,b,c}$ where $a\geq 0$ and
$1\leq b\leq c$ and $a+b+c=n-1$, or $G=\bigtriangleup_{p,q,r}$ where
$0\leq p\leq q\leq r$ and $p+q+r=n-3$.
\end{thm}

Let $\mathcal {G}(n)$ denote the class of simple graphs of order
$n$. Give a graph theoretic parameter $f(G)$ and a positive integer
$n$, the \emph{Nordhaus-Gaddum(\textbf{N-G}) Problem} is to
determine sharp bounds for: $(1)$ $f(G)+f(\overline{G})$ and $(2)$
$f(G)\cdot f(\overline{G})$, as $G$ ranges over the class $\mathcal
{G}(n)$, and characterize the extremal graphs. The Nordhaus-Gaddum
type relations have received wide investigations. Recently,
Aouchiche and Hansen published a survey paper on this subject, see
\cite{Aouchiche}.

Xu \cite{Xu} obtained the Nordhaus-Gaddum results for the Steiner
$2$-diameter of graphs. In Section $3$, we obtain the
Nordhaus-Gaddum results for the Steiner $k$-diameter of graphs.
\begin{thm}\label{th5}
Let $G\in \mathcal {G}(n)$ and let $k$ be an integer with $3\leq
k\leq n$. Then

$(i)$ $2k-1-x\leq sdiam_k(G)+sdiam_k(\bar{G})\leq
\max\{n+k-1,4k-2\}$;

$(ii)$ $(k-1)(k-x)\leq sdiam_k(G)\cdot sdiam_k(\bar{G})\leq
\max\{k(n-1),(2k-1)^2\}$,

where if $n\geq 2k-2$ then $x=0$; $x=1$ for any positive integer $n
$.
\end{thm}

For $k=n,n-1,n-2,3$, we improve the above Nordhaus-Gaddum results of
Steiner $k$-diameter and obtain the following results.
\begin{obs}\label{obs3}
Let $G$ be a graph of order $n \ (n\geq 3)$. Then

$(i)$ $sdiam_n(G)+sdiam_n(\overline{G})=2n-2$;

$(ii)$ $sdiam_n(G)\cdot sdiam_3(\overline{G})=(n-1)^2$.
\end{obs}

Akiyama and Harary \cite{AkiyamaH} characterized the graphs for
which $G$ and $\overline{G}$ both have connectivity one.
\begin{lem}{\upshape \cite{AkiyamaH}}\label{lem0}
Let $G$ be graph with $n$ vertices. Then
$\kappa(G)=\kappa(\bar{G})=1$ if and only if $G$ satisfies the
following conditions.

$(i)$ $\kappa(G)=1$ and $\Delta(G)=n-2$;

$(ii)$ $\kappa(G)=1$, $\Delta(G)\leq n-3$ and $G$ has a cutvertex
$v$ with pendant edge $e$ and pendant vertex $u$ such that $G-u$
contains a spanning complete bipartite subgraph.
\end{lem}

By Lemma \ref{lem0}, we obtain the following result.
\begin{pro}\label{proA}
Let $G$ be a graph of order $n \ (n\geq 5)$.

$(i)$ $2n-4\leq sdiam_{n-1}(G)+sdiam_{n-1}(\overline{G})\leq 2n-2$;

$(ii)$ $(n-2)^2\leq sdiam_{n-1}(G)\cdot
sdiam_{n-1}(\overline{G})\leq (n-1)^2$.

Moreover,

$(a)$ $sdiam_{n-1}(G)+sdiam_{n-1}(\overline{G})=2n-4$ or
$sdiam_{n-1}(G)\cdot sdiam_{n-1}(\overline{G})=(n-2)^2$ if and only
if both $G$ and $\overline{G}$ are $2$-connected;

$(b)$ $sdiam_{n-1}(G)+sdiam_{n-1}(\overline{G})=2n-3$ or
$sdiam_{n-1}(G)\cdot sdiam_{n-1}(\overline{G})=(n-1)(n-2)$ if and
only if $\lambda(G)=1$ and $\overline{G}$ are $2$-connected, or
$\lambda(\overline{G})=1$ and $G$ are $2$-connected.

$(c)$ $sdiam_{n-1}(G)+sdiam_{n-1}(\overline{G})=2n-2$ or
$sdiam_{n-1}(G)\cdot sdiam_{n-1}(\overline{G})=(n-1)^2$ if and only
if $G$ satisfies the following conditions.

$\bullet$ $\kappa(G)=1$, $\Delta(G)=n-2$;

$\bullet$ $\kappa(G)=1$, $\Delta(G)\leq n-3$ and $G$ has a cutvertex
$v$ with pendant edge $e$ and pendant vertex $u$ such that $G-u$
contains a spanning complete bipartite subgraph.
\end{pro}

\begin{pro}\label{proB}
Let $G$ be a graph of order $n \ (n\geq 5)$. If both $G$ and
$\overline{G}$ contains at least two cut vertices, then

$(i)$ $2n-6\leq sdiam_{n-2}(G)+sdiam_{n-2}(\overline{G})\leq 2n-2$;

$(ii)$ $(n-3)^2\leq sdiam_{n-2}(G)\cdot
sdiam_{n-2}(\overline{G})\leq (n-1)^2$.

Otherwise,

$(iii)$ $2n-6\leq sdiam_{n-2}(G)+sdiam_{n-2}(\overline{G})\leq
2n-3$;

$(iv)$ $(n-3)^2\leq sdiam_{n-2}(G)\cdot
sdiam_{n-2}(\overline{G})\leq (n-1)(n-2)$.

Moreover, the upper and lower bounds are sharp.
\end{pro}

For Steiner $3$-diameter, we improve the result in Theorem \ref{th5}
and prove the following result in Section 3.
\begin{pro}\label{proC}
Let $G$ be a graph of order $n \ (n\geq 10)$. Then

$(i)$ $6\leq sdiam_3(G)+sdiam_3(\overline{G})\leq n+2$;

$(ii)$ $9\leq sdiam_3(G)\cdot sdiam_3(\overline{G})\leq 3(n-1)$.

Moreover, the bounds are sharp.
\end{pro}

\section{Graphs with given Steiner $3$-diameter}

In this section, we characterize graphs with $sdiam_3(G)=2,3,n-1$
and give the proofs of Theorems \ref{th2}, \ref{th3} and \ref{th4}.

The following observation is easily seen.

\begin{obs}\label{obs3}
If $H$ is a spanning subgraph of $G$, then $sdiam_k(G)\leq
sdiam_k(H)$.
\end{obs}

When $G=T$ is a tree of order $n$, graphs attaining the upper bound
of Proposition \ref{pro1} can be characterized in the following,
which will be used later.

\begin{pro}\label{pro2}
Let $k,n$ be two integers with $2\leq k\leq n$, and let $T$ be a
tree of order $n$. Then $sdiam_k(T)=n-1$ if and only if $r\leq k$,
where $r$ is the number of leaves in $T$.
\end{pro}
\begin{pf}
Suppose $r\leq k$. Let $v_1,v_2,\cdots,v_{r}$ be all the leaves of
$T$. Choose $S\subseteq V(T)$ and $|S|=k$ such that
$v_1,v_2,\cdots,v_{r}\in S$. Then any $S$-Steiner tree must use all
edges of $T$. Since $|E(T)|=n-1$, it follows that $d_{T}(S)\geq
|E(T)|=n-1$ and hence $sdiam_k(T)\geq n-1$. Combining this with
Proposition \ref{pro1}, we have $sdiam_k(T)=n-1$.

Conversely, suppose $sdiam_k(T)=n-1$. If $s\geq k+1$, then for any
$S\subseteq V(G)$ with $|S|=k$, there exists a leaf $v$ in $T$ such
that $v\notin S$. Let $T'=T-v$. Then $T'$ is a $S$-Steiner tree and
hence $d_{T}(S)\leq n-2$. From the arbitrariness of $S$, we have
$sdiam_k(T)\leq n-2<n-1$, a contradiction. So $s\leq k$.\qed
\end{pf}

From Proposition \ref{pro1}, we have $k-1\leq sdiam_k(G)\leq n-1$.
We now show a property of the graphs attaining the lower bound.

\begin{lem}\label{lem1}
Let $n,k$ be two integers with $2\leq k\leq n$, and let $G$ be a
connected graph of order $n$. If $sdiam_k(G)=k-1$, then $0\leq
\Delta(\overline{G})\leq k-2$, namely, $n-k+1\leq \delta(G)\leq
n-1$.
\end{lem}
\begin{pf}
Suppose $\Delta(\overline{G})\geq k-1$. Then there exists a vertex
$u\in V(\overline{G})$ such that $d_{\overline{G}}(u)\geq k-1$. Pick
up $v_1,v_2,\cdots,v_{k-1}\in N_{\overline{G}}(u)$. Let
$S=\{u,v_1,v_2,\cdots,v_{k-1}\}$. Since $uv_i\in E(\overline{G}) \
(1\leq i\leq k-1)$, it follows that $uv_i\notin E(G)$ and hence $u$
is an isolated vertex in $G[S]$. Thus, any $S$-Steiner tree must use
$k$ edges of $E(G)$, which implies that $d_{G}(S)\geq k$. Therefore,
$sdiam_{k}(G)\geq k$, a contradiction. So $0\leq
\Delta(\overline{G})\leq k-2$, namely, $n-k+1\leq \delta(G)\leq
n-1$.\qed \vspace{3pt}
\end{pf}

\noindent{\bf Proof of Theorem \ref{th2}:} For Lemma \ref{lem1}, if
$sdiam_3(G)=2$, then $0\leq \Delta(\overline{G})\leq 1$. Conversely,
if $0\leq \Delta(\overline{G})\leq 1$, then $n-2\leq \delta(G)\leq
n-1$. Thus, $G$ is a graph obtained from the complete graph of order
$n$ by deleting some independent edges. For any
$S=\{u,v,w\}\subseteq V(G)$, at least two elements in $\{uv,vw,uw\}$
belong to $E(G)$. Without loss of generality, let $uv,vw\in E(G)$.
It is clear that the tree $T$ induced by the edges in $\{uv,vw\}$ is
an $S$-Steiner tree and hence $d_G(S)\leq 2$. From the arbitrariness
of $S$, we have $sdiam_3(G)\leq 2$ and hence $sdiam_3(G)=2$ by
Proposition \ref{pro1}. The proof is complete. \qed \vspace{3pt}

\noindent{\bf Proof of Theorem \ref{th3}:} Suppose that $G$ is a
graph with $sdiam_3(G)=3$. From Theorem \ref{th2}, we have
$\Delta(\overline{G})\geq 2$. It suffices to prove the following two
claims.

\textbf{Claim 1}. $\overline{G}$ does not contain a triple-star as
its subgraph.

Assume, to the contrary, that $\overline{G}$ contains a triple-star
$H_1$ as its subgraph. Choose $S=\{u,v,w\}$. Then $uv,uw,vw\in
E(\overline{G})$ and hence $uv,uw,vw\notin E(G)$. For any $x\in
V(G)-S$, one can see that $xu\notin E(G)$ or $xv\notin E(G)$ or
$xw\notin E(G)$. Observe that any $S$-Steiner tree $T$ must occupy
at least one vertex of $V(G)-S$, say $y$. Then $yu\notin E(G)$ or
$yv\notin E(G)$ or $yw\notin E(G)$. Without loss of generality, let
$yu\notin E(G)$. Therefore, the tree $T$ must occupy at least one
vertex of $V(G)-\{u,v,w,y\}$. Thus the tree $T$ contains at least
$5$ vertices in $G$, which implies that $d_{G}(S)\geq 4$ and hence
$sdiam_3(G)\geq 4$, a contradiction. So $\overline{G}$ does not
contain $H_1$ as its subgraph.

\textbf{Claim 2}. $\overline{G}$ does not contain $H_2$ as its
subgraph.

Assume, to the contrary, that $G$ contains $H_2$ as its subgraph.
Choose $S=\{u,v,w\}\subseteq V(G)$. Since $uv,vw\in
E(\overline{G})$, it follows that $uv,vw\notin E(G)$. Clearly, any
$S$-Steiner tree $T$ uses at least one vertex in $V(G)-S$. For each
$x\in V(G)-S$, we have $xu,xv,xw\in E(\overline{G})$ or $xu,xv\in
E(\overline{G})$ or $xv,xw\in E(\overline{G})$ or $xu,xw\in
E(\overline{G})$ or $xv\in E(\overline{G})$, that is,
$xu,xv,xw\notin E(G)$ or $xu,xv\notin E(G)$ or $xv,xw\notin E(G)$ or
$xu,xw\notin E(G)$ or $xv\notin E(G)$. One can see that the tree $T$
connecting $S$ uses at least two vertices in $V(G)-S$. Therefore,
$e(T)\geq 4$ and $d_{G}(S)\geq 4$, which results in $sdiam_3(G)\geq
4$, a contradiction. So $\overline{G}$ does not contain $H_2$ as its
subgraph.

From the above arguments, we know that the result holds.

Conversely, suppose that $G$ is a connected graph such that
$\Delta(\overline{G})\geq 2$ and $\overline{G}$ does not contain
both $H_1$ and $H_2$ as its subgraph. From the definition of
$sdiam_3(G)$, it suffices to show that $d_{G}(S)=3$ for any
$S\subseteq V(G)$. Set $S=\{u,v,w\}$. Then $0\leq |E(G[S])|\leq 3$.

If $2\leq |E(G[S])|\leq 3$, then there are two edges in $G[S]$
belonging to $E(G)$, say $uv,vw$. Therefore, the tree $T$ induced by
the edges in $\{uv,vw\}$ is an $S$-Steiner tree in $G$, which
results in $d_{G}(S)=2<3$, as desired.

Suppose $|E(G[S])|=0$. Then $uv,vw,uw\notin E(G)$ and hence
$uv,vw,uw\in E(\overline{G})$. Because $\overline{G}$ does not
contain the subgraph $H_1$ as its subgraph, there exists a vertex
$y\in V(G)-S$ such that $yu,yv,yw\notin E(\overline{G})$, which
implies $yu,yv,yw\in E(G)$. It is clear that the tree $T$ induced by
the edges in $\{yu,yv,yw\}$ is an $S$-Steiner tree in $G$ and hence
$d_{G}(S)\leq 3$, as desired.

Suppose $|E(G[S])|=1$. Without loss of generality, let $uw\in E(G)$.
Then $uv,vw\in E(\overline{G})$. Since $\overline{G}$ does not
contain $H_2$ as its subgraph, there exists a vertex $x\in V(G)-S$
such that $xu\in E(\overline{G})$ but $xv,xw\notin E(\overline{G})$,
or $xw\in E(\overline{G})$ but $xu,xv\notin E(\overline{G})$. By
symmetry, we only need to consider the former case. Then $xv,xw\in
E(G)$. Combining this with $uw\in E(G)$, the tree $T$ induced by the
edges in $\{xv,xw,uw\}$ is an $S$-Steiner tree in $G$, namely,
$d_{G}(S)\leq 3$, as desired.

From the arbitrariness of $S$, we know that $sdiam_3(G)\leq 3$.
Since $\Delta(\overline{G})\geq 2$, Theorem \ref{th2} implies that
$sdiam_3(G)=3$. The proof is now complete. \qed \vspace{3pt}

We are now in a position to give the proof of Theorem \ref{th4}.

\begin{lem}\label{lem2}
Let $G$ be a connected graph of order $n \ (n\geq 5)$. If $4\leq
c(G)\leq n$, then $sdiam_3(G)\leq n-2$, where $c(G)$ is the
circumference of the graph $G$.
\end{lem}
\begin{pf}
If $c(G)=n$, then there is a Hamilton cycle $C_n$ in $G$. From
Observations \ref{obs1} and \ref{obs3}, we have $sdiam_3(G)\leq
sdiam_3(C_n)=\lfloor\frac{2}{3}n\rfloor\leq n-2$. Let $c(G)=t \
(4\leq t\leq n-1)$. Then there exists a cycle of order $t$ in $G$,
say $C_t=v_1v_2\cdots v_tv_1$. Let $G_1,G_2,\cdots,G_r$ be the
connected components of $G-V(C_t)$.

Suppose $r\geq 4$. Clearly, each connected component $G_i \ (1\leq
i\leq r)$ contains a spanning tree $T_i$ (note that if $G_i$ is
trivial, then $T_i$ is trivial). Since $G$ is connected, there is an
edge $e_i$ such that one endpoint of $e_i$ belongs to $V(T_i)$ and
the other endpoint belongs to $V(C_t)$. Furthermore, we choose one
edge from the cycle $C_t$, say $e$, and delete it. Then the tree $T$
induced by the edges in $(\bigcup_{i=1}^r E(T_i))\cup
(\bigcup_{i=1}^r e_j)\cup (E(C_t)-e)$ is a spanning tree of $G$ with
at least four leaves. From Proposition \ref{pro2} and Observation
\ref{obs3}, $sdiam_3(G)\leq sdiam_3(T)\leq n-2$, as desired.

We now assume $r\leq 3$. It suffices to show that $d_G(S)\leq n-2$
for any $S\subseteq V(G)$ with $|S|=3$. We have the following four
cases to consider. If $|S\cap V(C_t)|=3$, then, from Observation
\ref{obs1}, $d_G(S)\leq sdiam_3(C_t)=\lfloor\frac{2}{3}t\rfloor\leq
\frac{2}{3}t\leq \frac{2}{3}(n-1)\leq n-2$, as desired. If $|S\cap
V(C_t)|=2$, then there exists a vertex $x\in S$ such that $x\in
V(G-V(C_t))$. Then $x$ must belong to some connected component in
$G-V(C_t)$. Without loss of generality, let $x\in V(G_1)$, and let
$S=\{x,v_i,v_j\}$ where $v_i,v_j\in V(C_t) \ (1\leq i\neq j\leq t)$.
Because $G_1$ is connected, $G_1$ contains a spanning tree, say
$T_1$. Since $G$ is connected, we can find an edge $e_1$ with one
endpoint belonging to $V(T_1)$ and the other, say $v_k$, belonging
to $V(C_t)$ (note that $v_k,v_i$ or $v_k,v_j$ are not necessarily
different). Since $d(\{v_i,v_j,v_k\})\leq
sdiam_3(C_t)=\lfloor\frac{2}{3}t\rfloor$, we have $d_G(S)\leq
d(\{v_i,v_j,v_k\})+|E(T_1)|+1=d(\{v_i,v_j,v_k\})+|V(T_1)|\leq
\lfloor\frac{2}{3}t\rfloor+n-t\leq n-\frac{1}{3}t<n-1$ and hence
$d_G(S)\leq n-2$, as desired.

Suppose $|S\cap V(C_t)|=1$. Then there exist two vertices $x,y\in S$
such that $x,y\in V(G-V(C_t))$. Set $S=\{x,y,v_i\}$ where $v_i\in
V(C_t) \ (1\leq i\leq t)$. Thus, $x,y$ must belong to the same
connected component of $G-V(C_t)$, or $x,y$ belong to two different
connected components. Consider the former case. Without loss of
generality, let $x,y\in V(G_1)$. Since $G_1$ is connected, it
follows that $G_1$ contains a spanning tree, say $T_1$. Because $G$
is connected, we can find an edge $e_1$ with one endpoint belonging
to $V(T_1)$ and the other, say $v_j$, belonging to $V(C_t)$ (note
that $v_i$ and $v_j$ are not necessarily different). Since
$d(\{v_i,v_j\})\leq \lfloor\frac{1}{2}t\rfloor$, it follows that
$d_G(S)\leq d(\{v_i,v_j\})+|E(T_1)|+1=d(\{v_i,v_j\})+|V(T_1)|\leq
n-t+\lfloor\frac{1}{2}t\rfloor\leq n-\lceil\frac{1}{2}t\rceil$.
Since $t\geq 4$, we have $d_G(S)\leq n-2$, as desired. Consider the
latter case. Without loss of generality, let $x\in V(G_1)$ and $y\in
V(G_2)$. Clearly, $G_i \ (i=1,2)$ contains a spanning tree $T_i$. We
can find the edges $e_1,e_2$ with one endpoint belonging to
$V(T_1),V(T_2)$ and the other, say $v_j,v_k$, belonging to $V(C_t)$,
respectively (note that $v_i,v_j,v_k$ are not necessarily
different). Since $d(\{v_i,v_j,v_k\})\leq
sdiam_3(C_t)=\lfloor\frac{2}{3}t\rfloor$, we have $d_G(S)\leq
d(\{v_i,v_j,v_k\})+|E(T_1)|+|E(T_2)|+2=
d(\{v_i,v_j,v_k\})+|V(T_1)|+|V(T_2)|\leq
\lfloor\frac{2}{3}t\rfloor+n-t\leq n-\lceil\frac{1}{3}t\rceil$ and
hence $d_G(S)\leq n-2$, as desired.

Suppose $|S\cap V(C_t)|=0$. Then $S\subseteq V(G-V(C_t))$. Let
$S=\{x,y,z\}$. Thus, $x,y,z$ belong to three different connected
components, or $x,y,z$ belong to two different connected components,
or $x,y,z$ must belong to one connected component. We only prove the
first case, the other two cases can be proved similarly. Without
loss of generality, let $x\in V(G_1)$, $y\in V(G_2)$ and $z\in
V(G_3)$. For $i=1,2,3$, $G_i$ contains a spanning tree $T_i$. Since
$G$ is connected, we can find the edges $e_1,e_2,e_3$ with one
endpoint belonging to $V(T_1),V(T_2),V(T_3)$ and the other, say
$v_i,v_j,v_k$, belonging to $V(C_t)$, respectively (note that
$v_j,v_k,v_j$ are not necessarily different). Since
$d_G(\{v_i,v_j,v_k\})\leq sdiam_3(C_t)=\lfloor\frac{2}{3}t\rfloor$,
we have $d_G(S)\leq d(\{v_i,v_j,v_k\})+\sum_{i=1}^3|E(T_i)|+3=
d(\{v_i,v_j,v_k\})+\sum_{i=1}^3|V(T_i)|\leq
\lfloor\frac{2}{3}t\rfloor+n-t\leq n-\lceil\frac{1}{3}t\rceil$ and
hence $d_G(S)\leq n-2$, as desired.

From the above arguments, we conclude that $sdiam_3(G)\leq n-2$. The
proof is now complete. \qed
\end{pf}\vspace{3pt}

If $T$ is a nontrivial tree and $S\subseteq V(T)$, where $|S|\geq
2$, then there is a unique subtree $T_s$ of size $d(S)$ containing
the vertices of $S$. We refer to such a tree as the \emph{tree
generated by $S$}.

Chartrand, Oellermann, Tian and Zou \cite{Chartrand} obtained the
following result.
\begin{lem}{\upshape \cite{Chartrand}}\label{lemF}
If $H$ is a subgraph of a graph $G$ and $v$ is a vertex of $G$, then
$d(v,H)$ denotes the minimum distance from $v$ to a vertex of $H$.
Therefore,
$$
d(S\cup \{v\})=d(S)+d(v,T_s).
$$
\end{lem}

\noindent{\bf Proof of Theorem \ref{th4}:} For $n=3$,
$sdiam_3(G)=n-1=2$ if and only if $G=P_3=T_{0,1,1}$ or
$G=K_3=\bigtriangleup_{0,0,0}$. For $n=4$, $sdiam_3(G)=n-1=3$ if and
only if $G=P_4=T_{0,1,2}$ or $G=\bigtriangleup_{0,0,1}$. We now
assume $n\geq 5$.

Suppose $G=T_{a,b,c}$ where $0\leq a\leq b\leq c$ and $1\leq b\leq
c$ and $a+b+c=n-1$. Since there are at most three leaves in $G$, it
follows from Proposition \ref{pro2} that $sdiam_3(G)=n-1$. Suppose
$G=\bigtriangleup_{p,q,r}$ where $0\leq p\leq q\leq r$ and
$p+q+r=n-3$. From Proposition \ref{pro1}, we have $sdiam_3(G)\leq
n-1$. It suffices to show that $sdiam_3(G)\geq n-1$. Choose the
three leaves in $T_{a,b,c}$, say $x,y,z$, such that $x\in V(P_a)$,
$y\in V(P_b)$ and $z\in V(P_c)$. Let $S'=\{x,z\}$ and $S=\{x,y,z\}$.
From Lemma \ref{lemF}, $d_G(S)=d_G(S'\cup
\{y\})=d_G(S')+d(y,T_s)=(b+c)+a=n-1$, and hence
$sdiam_3(G)=sdiam_3(T_{a,b,c})=n-1$, as desired. Similarly, we can
get that $sdiam_3(\bigtriangleup_{p,q,r})=n-1$, as desired.

Conversely, suppose $sdiam_3(G)=n-1$. If $G$ is a tree, then, by
Proposition \ref{pro2}, $G$ contains at most three leaves. Thus,
$G=T_{a,b,c}$, where $0\leq a\leq b\leq c$ and $1\leq b\leq c$ and
$a+b+c=n-1$. Now, we consider the graph $G$ containing cycles.
Recall that $c(G)$ is the circumference of the graph $G$. Obviously,
$3\leq c(G)\leq n$. If $4\leq c(G)\leq n$, then it follows from
Lemma \ref{lem2} that $sdiam_3(G)\leq n-2$, a contradiction.
Therefore, $c(G)=3$. Suppose that $G$ contains at least two
triangles. If there exist two triangles having at most one common
vertex, then $G$ contains a spanning tree with at least four leaves,
say $T$. From Observation \ref{obs3} and Proposition \ref{pro1}, we
have $sdiam_3(G)\leq sdiam_3(T)\leq n-2$, a contradiction. So we
assume that there exist two triangles having two common vertices in
$G$. Therefore, $G$ contains $K_4^-$ as its subgraph, where $K_4^-$
is a graph obtained from a clique $K_4$ by deleting one edge. Now,
we consider the two vertices of degree $3$ in $K_4^-$. If the degree
of each such vertex in $K_4^-$ is larger than $4$ in $G$, then $G$
contains a spanning tree with four leaves. Again from Observation
\ref{obs3} and Proposition \ref{pro1}, $sdiam_3(G)\leq
sdiam_3(T)\leq n-2$, a contradiction. Then $G$ contains the graph
$H$ as its subgraph, where $H$ is a graph obtained from $K_4^-$ and
two paths by identifying one endvertex of each path and each vertex
of degree $2$ in $K_4^-$. One can also that $sdiam_3(H)\leq n-2$ and
hence $sdiam_3(G)\leq sdiam_3(H)\leq n-2$, a contradiction. From the
above arguments, we conclude that $G$ only contains one triangle and
hence $G=\bigtriangleup_{a,b,c}$. The proof is complete.\qed

\section{Nordhaus-Gaddum results}

The following proposition is a preparation of the proof of Theorem
\ref{th5}.
\begin{pro}\label{pro6}
Let $G$ be a connected graph. If $sdiam_k(G)\geq 2k$, then
$sdiam_k(\overline{G})\leq k$.
\end{pro}
\begin{pf}
For any $S\subseteq V(G)$ and $|S|=k$, if $G[S]$ is not connected,
then $\overline{G}[S]$ is connected, and hence
$d_{\overline{G}}(S)=k-1<k$. Suppose that $G[S]$ is connected. Then
we have the following claim.

\textbf{Claim 1.} There exists a vertex $u\in V(G)-S$ such that
$|E_G[u,S]|=0$.

Assume, to the contrary, that $|E_G[x,S]|\geq 1$ for any $x\in
V(G)-S$. For any $S'\subseteq V(G)$ and $|S'|=k$, since $G[S]$ is
connected and $|E_G[x,S]|\geq 1$ for any $x\in S'-S$, it follows
that $G[S\cup S']$ is connected, and hence $d_{G}(S')\leq 2k-1$.
From the arbitrariness of $S'$, we have $sdiam_k(G)\leq 2k-1$, a
contradiction. \qed

From Claim 1, there exists a vertex $u\in V(G)-S$ such that
$|E_{\overline{G}}[u,S]|=k$, and the tree induced by these $k$ edges
is an $S$-Steiner tree in $\overline{G}$. So
$d_{\overline{G}}(S)=k$. From the arbitrariness of $S$, we have
$sdiam_k(\overline{G})\leq k$, as desired.\qed
\end{pf}\vskip 0.2cm

\noindent{\bf Proof of Theorem \ref{th5}:} We first give the proof
of the upper bounds. If $sdiam_k(G)\geq 2k$, then it follows from
Proposition \ref{pro6} that $sdiam_k(\overline{G})\leq k$.
Furthermore, since $sdiam_k(G)\leq n-1$, we have
$sdiam_k(G)+sdiam_k(\overline{G})\leq n+k-1$ and $sdiam_k(G)\cdot
sdiam_k(\overline{G})\leq k(n-1)$. By the same reason, if
$sdiam_k(\overline{G})\geq 2k$, then $sdiam_k(G)\leq k$, and hence
$sdiam_k(G)+sdiam_k(\overline{G})\leq n+k-1$ and $sdiam_k(G)\cdot
sdiam_k(\overline{G})\leq k(n-1)$. We now assume that
$sdiam_k(G)\leq 2k-1$ and $sdiam_k(\overline{G})\leq 2k-1$. Then
$sdiam_k(G)+sdiam_k(\overline{G})\leq 4k-2$, and hence
$sdiam_k(G)+sdiam_k(\overline{G})\leq \max\{n+k-1,4k-2\}$ and
$sdiam_k(G)\cdot sdiam_k(\overline{G})\leq \max\{k(n-1),(2k-1)^2\}$.

Next, we show the proof of the lower bounds. From Proposition
\ref{pro1}, since $sdiam_k(G)\geq k-1$ and
$sdiam_k(\overline{G})\geq k-1$, we have
$sdiam_k(G)+sdiam_k(\overline{G})\geq 2k-2$ and $sdiam_k(G)\cdot
sdiam_k(\overline{G})\geq (k-1)^2$. Since $n\geq 2k-2$, we claim
that $sdiam_k(G)+sdiam_k(\overline{G})\geq 2k-1$ and
$sdiam_k(G)\cdot sdiam_k(\overline{G})\geq (k-1)k$. Assume, to the
contrary, that $sdiam_k(G)+sdiam_k(\overline{G})=2k-2$ and
$sdiam_k(G)\cdot sdiam_k(\overline{G})\geq (k-1)^2$. Then
$sdiam_k(G)=sdiam_k(\overline{G})=k-1$. From Lemma \ref{lem1}, we
have $n-k+1\leq \delta(G)\leq n-1$ and $0\leq \Delta(G)\leq k-2$,
and hence $n\leq 2k-3$, a contradiction. So
$sdiam_k(G)+sdiam_k(\overline{G})\geq 2k-1$ and $sdiam_k(G)\cdot
sdiam_k(\overline{G})\geq (k-1)k$.\qed\vskip 0.2cm

\begin{lem}\label{lemM}
Let $G$ be a graph. Then $sdiam_{n-1}(G)=n-2$ if and only if $G$ is
$2$-connected.
\end{lem}
\begin{pf}
Suppose that $G$ is $2$-connected. For any $S\subseteq V(G)$ and
$|S|=n-1$, there exists a unique vertex $V(G)-S$, say $v$, such that
$G-v$ is connected, and hence $G-v$ contains a spanning tree, which
implies $d_G(S)\leq n-2$. From the arbitrariness of $S$, we have
$sdiam_{n-1}(G)\leq n-2$. From Proposition \ref{pro1},
$sdiam_{n-1}(G)=n-2$.

Conversely, we suppose $sdiam_{n-1}(G)=n-2$. If $G$ is not
$2$-connected, then there exists a cut vertex in $G$, say $v$.
Choose $S=V(G)-v$. Then $|S|=n-1$. Observe that any $S$ Steiner tree
must use all the vertices of $G$. Thus $d_G(S)\geq n-1$, which
contradicts $sdiam_{n-1}(G)=n-2$.\qed
\end{pf}\vskip 0.2cm

By Proposition \ref{pro6} and Lemma \ref{lemM}, we can give the
proof of Proposition \ref{proA}.

\noindent{\bf Proof of Proposition \ref{proA}:} From Proposition
\ref{pro1}, we have $2n-4\leq
sdiam_{n-1}(G)+sdiam_{n-1}(\overline{G})\leq 2n-2$ and $(n-2)^2\leq
sdiam_{n-1}(G)\cdot sdiam_{n-1}(\overline{G})\leq (n-1)^2$. Clearly,
$sdiam_{n-1}(G)+sdiam_{n-1}(\overline{G})=2n-4$ or
$sdiam_{n-1}(G)\cdot sdiam_{n-1}(\overline{G})=(n-2)^2$ if and only
if $sdiam_{n-1}(G)=sdiam_{n-1}(\overline{G})=n-2$. From Lemma
\ref{lemM}, $sdiam_{n-1}(G)+sdiam_{n-1}(\overline{G})=2n-4$ or
$sdiam_{n-1}(G)\cdot sdiam_{n-1}(\overline{G})=(n-2)^2$ if and only
if both $G$ and $\overline{G}$ are $2$-connected.

It is clear that $sdiam_{n-1}(G)+sdiam_{n-1}(\overline{G})=2n-3$ or
$sdiam_{n-1}(G)\cdot sdiam_{n-1}(\overline{G})=(n-1)(n-2)$ if and
only if $sdiam_{n-1}(G)=n-2$ and $sdiam_{n-1}(\overline{G})=n-1$, or
$sdiam_{n-1}(G)=n-1$ and $sdiam_{n-1}(\overline{G})=n-2$.
Furthermore, $sdiam_{n-1}(G)+sdiam_{n-1}(\overline{G})$ $=2n-3$ or
$sdiam_{n-1}(G)\cdot sdiam_{n-1}(\overline{G})=(n-1)(n-2)$ if and
only if $\lambda(G)=1$ and $\overline{G}$ are $2$-connected, or
$\lambda(\overline{G})=1$ and $G$ are $2$-connected.

For the remaining case,
$sdiam_{n-1}(G)+sdiam_{n-1}(\overline{G})=2n-2$ or
$sdiam_{n-1}(G)\cdot sdiam_{n-1}(\overline{G})=(n-1)^2$ if and only
if $sdiam_{n-1}(G)=sdiam_{n-1}(\overline{G})=n-1$. From Lemma
\ref{lem0}, $sdiam_{n-1}(G)+sdiam_{n-1}(\overline{G})=2n-2$ or
$sdiam_{n-1}(G)\cdot sdiam_{n-1}(\overline{G})=(n-1)^2$ if and only
if $G$ satisfies the following conditions.

$\bullet$ $\kappa(G)=1$, $\Delta(G)=n-2$;

$\bullet$ $\kappa(G)=1$, $\Delta(G)\leq n-3$ and $G$ has a cutvertex
$v$ with pendant edge $e$ and pendant vertex $u$ such that $G-u$
contains a spanning complete bipartite subgraph.\qed\vskip 0.2cm

\noindent{\bf Proof of Proposition \ref{proB}:} From Proposition
\ref{pro1}, $2n-6\leq sdiam_{n-2}(G)+sdiam_{n-2}(\overline{G})\leq
2n-2$ and $(n-3)^2\leq sdiam_{n-2}(G)\cdot
sdiam_{n-2}(\overline{G})\leq (n-1)^2$. So the results follow for
the case that both $G$ and $\overline{G}$ contain at least two cut
vertices. From now on, we assume that $G$ or $\overline{G}$ contains
only one cut vertex, or $G$ or $\overline{G}$ is $2$-connected.
Without loss of generality, we assume that $G$ contains only one cut
vertex or $G$ is $2$-connected. For any $S\subseteq V(G)$ and
$|S|=n-2$, there exists a vertex $v\in V(G)-S$ such that $G-v$ is
connected, and hence $G-v$ contains a spanning tree, which implies
$d_G(S)\leq n-2$. From the arbitrariness of $S$, we have
$sdiam_{n-2}(G)\leq n-2$. From Proposition \ref{pro1}, we have
$sdiam_{n-2}(\overline{G})\leq n-1$. So
$sdiam_{n-2}(G)+sdiam_{n-2}(\overline{G})\leq 2n-3$ and
$sdiam_{n-2}(G)\cdot sdiam_{n-2}(\overline{G})\leq (n-1)(n-2)$.
\qed\vskip 0.2cm

To show the sharpness of the bounds in Proposition \ref{proB}, we
consider the following example.

\noindent \textbf{Example 1:} Let $G=P_4$. Then $\overline{G}=P_4$,
$sdiam_{2}(P_4)=sdiam_{2}(\overline{P_4})=3$. Therefore, we have
$sdiam_{2}(P_4)+sdiam_{2}(\overline{P_4})=6=2n-4$ and
$sdiam_{2}(P_4)\cdot sdiam_{2}(\overline{P_4})=9=(n-1)^2$, which
implies that the upper bounds are sharp for the case both $G$ and
$\overline{G}$ contain at least two cut vertices. Let $S^*$ be a
tree obtained from a star of order $n-2$ and a path of length $2$ by
identifying the center of the star and a vertex of degree one in the
path. Then $\overline{S^*}$ is a graph obtained from a clique of
order $n-1$ by deleting an edge $uv$ and then adding an pendent edge
$vw$ at $v$. Choose $S=V(G)-\{u,w\}$. Then any $S$-Steiner tree uses
all the vertices of $V(G)$, and hence $d_G(S)\geq n-1$. From the
arbitrariness of $S$, we have $sdiam_{n-2}(G)\geq n-1$, and hence
$sdiam_{n-2}(G)=n-1$ by Proposition \ref{pro1}. Choose $S\subseteq
V(\overline{G})-w$ and $|S|=n-2$. Then any $S$-Steiner tree uses
$n-1$ vertices of $V(\overline{G})$, and hence
$d_{\overline{G}}(S)\geq n-2$. From the arbitrariness of $S$, we
have $sdiam_{n-2}(\overline{G})\geq n-2$. One can easily check that
$sdiam_{n-2}(\overline{G})\leq n-2$. So
$sdiam_{n-2}(\overline{G})=n-2$, and hence
$sdiam_{n-2}(G)+sdiam_{n-2}(\overline{G})=2n-3$ and
$sdiam_{n-2}(G)\cdot sdiam_{n-2}(\overline{G})=(n-1)(n-2)$. This
implies that the upper bounds in Proposition \ref{proB} are sharp.
Let $G$ be a graph such that both $G$ and $\overline{G}$ are
$3$-connected. For any $S\subseteq V(G)$ and $|S|=n-2$, there exist
two vertices $u,v$ in $V(G)-S$ such that $G-\{u,v\}$ is connected,
and hence $G-\{u,v\}$ contains a spanning tree, which implies
$d_G(S)\leq n-3$. From the arbitrariness of $S$, we have
$sdiam_{n-2}(G)\leq n-3$, and hence $sdiam_{n-2}(G)=n-3$ by
Proposition \ref{pro1}. Similarly, we have
$sdiam_{n-2}(\overline{G})=n-3$. Then
$sdiam_{n-2}(G)+sdiam_{n-2}(\overline{G})=2n-6$ and
$sdiam_{n-2}(G)\cdot sdiam_{n-2}(\overline{G})=(n-3)^2$, which
implies that the lower bounds in Proposition \ref{proB} are
sharp.\vskip 0.2cm

\noindent{\bf Proof of Proposition \ref{proC}:} The upper bounds
follow from Theorem \ref{th5}. We now show the lower bounds of
$sdiam_3(G)+sdiam_3(\overline{G})$ and $sdiam_3(G)\cdot
sdiam_3(\overline{G})$. If $sdiam_3(G)+sdiam_3(\overline{G})<6$ or
$sdiam_3(G)\cdot sdiam_3(\overline{G})<9$, then we have
$sdiam_3(G)=2$ or $sdiam_3(\overline{G})=2$. Without loss of
generality, let $sdiam_3(G)=2$. From Theorem \ref{th2}, we have
$0\leq \Delta(\overline{G})\leq 1$ and hence $\overline{G}$ is
disconnected. Thus $sdiam_3(\overline{G})=\infty$, which results in
$sdiam_3(G)+sdiam_3(\overline{G})=\infty$ and $sdiam_3(G)\cdot
sdiam_3(\overline{G})=\infty$, a contradiction. So
$sdiam_3(G)+sdiam_3(\overline{G})\geq 6$ and $sdiam_3(G)\cdot
sdiam_3(\overline{G})\geq 9$.\qed\vskip 0.2cm

To show the sharpness of the bounds in Proposition \ref{proC}, we
consider the following example.

\noindent \textbf{Example 2:} One can check that $G=P_n$ is a sharp
example for the upper bounds of this theorem. To show the sharpness
of the lower bounds, we consider the following example. If
$sdiam_3(G)+sdiam_3(\overline{G})=6$, then
$sdiam_3(G)=sdiam_3(\overline{G})=3$. Let $G'$ be a graph of order
$n-4$, and let $a,b,c,d$ be a path. Let $G$ be the graph obtained
from $G'$ and the path by adding edges between the vertex $a$ and
all vertices of $G'$ and adding edges between the vertex $d$ and all
vertices of $G'$; see Figure 4 $(a)$. We now show that
$sdiam_3(G)=sdiam_3(\overline{G})=3$. Choose $S=\{a,b,d\}$. Then it
is easy to see that $d_{G}(S)\geq 3$ and hence $sdiam_3(G)\geq 3$.
It suffices to prove that $d_{G}(S)\leq 3$ for any $S\subseteq V(G)$
with $|S|=3$. Suppose $|S\cap V(G')|=3$. Without loss of generality,
let $S=\{x,y,z\}$. Then the tree $T$ induced by the edges in $\{xa,
ya,za\}$ is an $S$-Steiner tree and hence $d_{G}(S)\leq 3$. Suppose
$|S\cap V(G')|=2$. Without loss of generality, let $x,y\in S\cap
V(G')$. If $a\in S$, then the tree $T$ induced by the edges in
$\{xa,ya\}$ is an $S$-Steiner tree, which implies $d_{G}(S)\leq 2$.
If $b\in S$, then the tree $T$ induced by the edges in
$\{xa,ya,ab\}$ is an $S$-Steiner tree and hence $d_{G}(S)\leq 3$.
\begin{figure}[!hbpt]
\begin{center}
\includegraphics[scale=0.8]{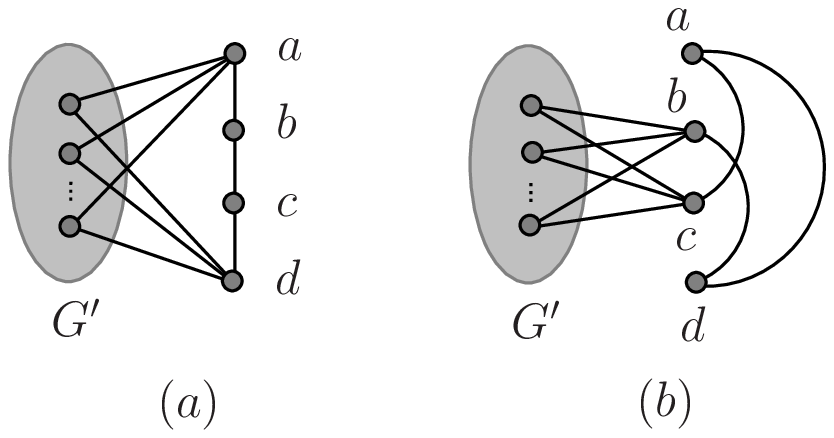}
\end{center}
\begin{center}
Figure 4: Graphs for Theorem \ref{th5}.
\end{center}\label{fig1}
\end{figure}
Suppose $|S\cap V(G')|=1$. Without loss of generality, let $x\in
S\cap V(G')$. If $a,b\in S$, then the tree $T$ induced by the edges
in $\{xa,ab\}$ is an $S$-Steiner tree and hence $d_{G}(S)\leq 2$. If
$b,c\in S$, then the tree $T$ induced by the edges in $\{xd,cd,bc\}$
is an $S$-Steiner tree and hence $d_{G}(S)\leq 3$. If $a,c\in S$,
then the tree $T$ induced by the edges in $\{xa,ab,bc\}$ is an
$S$-Steiner tree, which implies $d_{G}(S)\leq 3$. Suppose $|S\cap
V(G')|=0$. If $a,b,c\in S$, then the tree $T$ induced by the edges
in $\{ab,bc\}$ is an $S$-Steiner tree and hence $d_{G}(S)\leq 2$. If
$a,b,d\in S$, then the tree $T$ induced by the edges in
$\{ab,bc,cd\}$ is an $S$-Steiner tree, which implies $d_{G}(S)\leq
3$. From the arbitrariness of $S$, we conclude that $sdiam_3(G)\leq
3$ and hence $sdiam_3(G)=3$. Similarly, one can also check that
$sdiam_3(\overline{G})=3$. \qed

\end{document}